\documentclass[12pt, leqno, final]{amsart}

\usepackage{amssymb}
\usepackage{xypic}
\usepackage{showkeys}

\newcommand{\secnumber}{\S\arabic{section}}
\newcommand{\msection}[1]{\section{\bf{\hspace{-.43cm}\secnumber\ \  #1}}}

\swapnumbers
\newtheorem{theorem}{Theorem}[section]
\newtheorem{lemma}[theorem]{Lemma}
\newtheorem{proposition}[theorem]{Proposition}

\theoremstyle{definition}

\newtheorem{corollary}[theorem]{Corollary}

\newcommand{\symm}[1][d]{C^{(#1)}}
\newcommand{\w}[1][d]{W^{#1}}
\newcommand{\smooth}{{\rm smooth}}
\newcommand{\Gauss}{{\mathcal G}}

\newcommand{\Pdual}{({\mathbb P}^{g-1})^*}

\newcommand{\grass}[2]{{\rm G}(#1,#2)}
\newcommand{\f}[2]{{\rm F}({#1},{#2})}
\newcommand{\e}[2]{{\rm E}({#1},{#2})}

\newcommand{\lin}[1]{{\overline{#1}}}
\newcommand{\lint}[2]{{\overline{{#1}\cup{#2}}}}
\newcommand{\lind}[1]{{\lin{\phi_K(#1)}}}
\newcommand{\lintd}[2]{{\lint{\phi_K(#1)}{\phi_K(#2)}}}
\newcommand{\ab}[1]{{u(#1)}}
\newcommand{\abt}[2]{{(\ab{#1},\ab{#2})}} 
\newcommand{\pde}[2]{{\frac{\partial #1}{\partial #2}}}
\newcommand{\Grass}[1]{{\mathbb G}(#1)}
\newcommand{\m}{{\rm main}}

\begin{document}

\title{ A Generalized Torelli Theorem}

\author{Ajneet Dhillon}

\begin{abstract}{Given a smooth projective curve $C$ of positive
genus $g$, Torelli's
theorem asserts that the pair $(J(C),W^{g-1})$ determines
$C$. We show that the theorem is true with $W^{g-1}$ replaced
by $W^d$ for each $d$, in the range $1\le d\le g-1$. }
\end{abstract}

\maketitle

\msection{Introduction}

All curves in subsequent sections will be assumed to be smooth
projective curves over  ${\mathbb C}$. The genus of $C$ will
always be denoted by $g$. If $C$ is such a curve (with $g>0$)
we will let $J(C)$ denote its Jacobian and 
\[
u\colon C\rightarrow J(C)
\]
will be the Abel-Jacobi map. We will let $\symm$
denote the $d$th symmetric power of $C$ and for 
$1\le d\le g-1$, $\w$ will be the image of $\symm$ inside
the Jacobian under the Abel-Jacobi map. 
Since by a theorem of Riemann, the theta
divisor is a translate of $\w[g-1]$, Torelli's theorem asserts that the
pair $(J(C),\w[g-1])$ determines the curve, meaning that
if $C'$ is another curve such that there is an 
isomorphism $J(C)\cong J(C')$ carrying theta divisors
to theta divisors then the curves must be isomorphic.
Our aim is to show that an analogous statement holds
for each $1\le d<g-1$. With this in mind we will assume
in all following sections
 that $g\ge 4$, as smaller genera
are covered by existing theorems. Our strategy is largely 
based on the strategy in \cite{andre}.

As a corollary we have that two curves are isomorphic 
if and only if their $d$th symmeteric powers are isomorphic,
where $d$ is an integer smaller than the genus of one (and hence both)
of the curves.

This problem and the above mentioned Corollary
was originally proposed by Prof.~Donu Arapura.
Thanks also to the particpants of the Working Algebraic 
Geometry Seminar at Purdue, in particular to Prof. Kenji
Matsuki who pointed out a mistake in an earlier version.

\msection{Preliminaries}
\label{S:prelim}
The Jacobian of a curve $C$ is defined to be 
\[ J(C)={\rm H}^0(C,\Omega^1_C)^*/{\rm H}_1(C,\mathbb Z). \]
The Abel-Jacobi map is defined by
\[ \begin{array}{cc} 
	u:& C\longrightarrow J(C) \\
	  & p\longmapsto \int_{p_0}^p
   \end{array}
\]
where $p_0$ is a fixed basepoint. Let $\symm=C^d/S^d$ be the $d$th 
symmetric power of $C$. 
We identify the points of $\symm$ with 
effective divisors of degree $d$ on $C$.
 The Abel-Jacobi map
can be extended to a morphism
\[ u:\symm\longrightarrow J(C). \]
We have
\begin{theorem}[Abel's]
\label{T:abels}
	Let $D,D^\prime\in\symm$. Then 
\[
	 D\sim D^\prime\quad if\ and\ only\ if\quad u(D)=u(D^\prime)
\]
where the relation $\sim$ is linear equivalence.
\end{theorem}
\begin{proof}
See \cite{griffithsharris}.
\end{proof}
We let $\w=u(\symm)$. By Abel's Theorem $\w$ parameterises complete
linear systems of degree $d$ on $C$. Our aim is to reconstruct $C$ from
the pair $(J(C),\w)$ where $0<d\le g-1$. The main tool in doing this
will be the Gauss map, defined as follows. Take  $p\in\w_\smooth$
and let ${\rm T}_p(\w)$ be its holomorphic tangent space. There is an
automorphism, translation by $-p$, 
\[ 
   \begin{array}{cc}
	\tau_p:&J(C)\longrightarrow J(C) \\
	       &x\longmapsto x-p
   \end{array}
\]
This allows us to canonically identify ${\rm T}_p(\w)$ with a $d$-dimensional
subspace of ${\rm T}_0(J(C))\simeq {\rm H}^0(C,\Omega^1_C)^*$. This defines
the Gauss map
\[
  \Gauss:\w_\smooth\longrightarrow{\mathbb G}(d-1,g-1),
\]
where $\Grass{d-1,g-1}$ is the Grassmanian parameterizing $d-1$ dimensional
linear subvarieties  of ${\mathbb P}^{g-1}$ (or equivalently $d$-dimensional
subspaces of ${\mathbb C}^g$. The result we need is:

\begin{theorem}
\label{T:Gaussmap}
Let $\phi_K\colon C\rightarrow\Pdual$ be the canonical
morphism and let $D\in\symm$. Then $u(D)\in\w_\smooth$ if and only if 
${\rm dim}|D|=0$. If we denote by $\overline{\phi_K(D)}$ the linear
span of $D$ on the canonical curve then
\[
  \Gauss(u(D))=\overline{\phi_K(D)}.
\]
\end{theorem}
\begin{proof}
This result can be found in \S 2.7 of \cite{griffithsharris}.
\end{proof}

Note that the linear span of a multiple of a point is the appropriate
osculating plane to $C$ inside ${\mathbb P}^{g-1}$. The condition that
${\rm dim}|D|=0$ forces $\overline{\phi_K(D)}$ to be a $d-1$ dimensional
linear subvariety of ${\mathbb P}^{g-1}$. This is by

\begin{theorem}[Geometric Riemann-Roch]
\label{T:geometricRR}
For $D$ as in the above discussion we have
${\rm dim}|D|=d-1-{\rm dim}\overline{\phi_K(D)}.$
\end{theorem}	
\begin{proof}
Again this can be found in \cite{griffithsharris}.
\end{proof}

\msection{Our Strategy}

We first describe the idea behind the proof of the Torelli
theorem for curves, due to A. Andreotti, see \cite{andre}.
The Gauss map
\[
\Gauss\colon\w[g-1]_\smooth\rightarrow\Pdual
\]
is a quasi-finite morphism of degree
\[
{2g-2 \choose g-1}.
\]
To see this, a  hyperplane $H$ intersects the
image of a curve $C$ under its canonical morphism in
$2g-2$ points $p_1,p_2,\ldots,p_{2g-2}$, which are in
general posiion for a generic $H$.
By Theorem
\ref{T:Gaussmap} the fibre over $H$ consists of all images of divisors of
the form $u(p_{i_1}+p_{i_2}+\ldots +p_{i_{g-1}})$
where $i_j$ range over $\{1,2,\ldots,2g-2\}$.
If $C$ is non-hyperelliptic then
let $C^*$ be the dual variety to $C$, that is the locus
of all tangent hyperplanes to $\phi_K(C)$ inside
$\Pdual$. Now one would expect that the (closure of the) branch locus 
of $\Gauss$ to be $C^*$ since the fibre over a tangent hyperplane $H$
should have cardinality smaller than 
\[
{2g-2 \choose g-1}.
\]
(Since $H.C=2p_1+\ldots p_{2g-3}$, the first point is repeated
and there are fewer choices for points in the fibre.) It is known 
how to recover $C$ from $C^*$, for example see \cite{harris}.
In the case that $C$ is hyperelliptic the canonical morphism
$\phi_K\colon C\rightarrow{\mathbb P}^{g-1}$ is branched at
$2g+2$ points labelled $b_1,\ldots,b_{2g+2}$. We denote
by $C^*$ the dual variety to the rational normal curve
$\phi_K(C)$ and $b_i^*$ denotes the locus of all hyperplanes
passing through $b_i$. In the hyperelliptic case, by the same
reasoning as in the non-hyperelliptic case, one would 
expect that the branch locus of $\Gauss$ to be  
$C^*\cup b_1^*\cup\ldots b_{2g+2}^*$. It is known 
how to recover $C$ from this information. 

We would like to try to apply this technique to our situation.
Firstly, we may reduce to the case where $(g-1)/2<d<g-1$. 
To do this choose an integer $n$ so that $(g-1)/2<nd\le g-1$.
Then 
\[
W^{nd}=\underbrace{\w+\w+\ldots+\w}_{n\ {\rm times}}.
\]
The above addition is addition inside the Jacobian. 

Fix ${\mathbb P}^{g-1}={\mathbb P}({\rm H}^0(C,\Omega^1_C)^*)$.
Now consider the locus
\[
\begin{array}{c}
\f{d}{g}=\{(V,W)\in\Grass{d-1,{\mathbb P}^{g-1}}\times
\Grass{d-1,{\mathbb P}^{g-1}}\mid \\ \overline{V + W}\ne
{\mathbb P}^{g-1} \}.
\end{array}
\]
The notation $\overline{V + W}$ means linear span of $V$ and $W$.
So $\f{d}{g}$ is the locus of all pairs of $(d-1)$-dimensional 
linear subvarieties that are contained inside some hyperplane.
There is a rational morphism 
\[
\xymatrix{
\alpha\colon\f{d}{g} \ar@{-->}[r] & \Pdual
}
\]
defined by $(V,W)\mapsto \overline{V+W}$. 
We take $\e{d}{g}$ to be the pullback
of $\f{d}{g}$ under 
\[
\Gauss\times\Gauss\colon \w_\smooth\times
\w_\smooth\rightarrow\Grass{d-1,g-1}\times\Grass{d-1,g-1}.
\]
Now let $\beta$ be the composed rational morphism
\[ 
\xymatrix{
\beta\colon \e{d}{g} \ar@{-->}[r] & \Pdual.
}
\]
Arguing as in the case $d=g-1$ we see that the
branch locus of $\beta$ contains enough information
to recover $C$. Note that the hypothesis
$(g-1)/2<d<g-1$ is required to insure that $\e{d}{g}$
is not empty.

\msection{Generic Determinental Varieties}

Two identities that will be useful later are presented in this
section.

In this section 
 $d$ and $g$ will be non-negative integers with $(g-1)/2<d<g-1$.
We will need the case $g\ge 4$ later. 
Let $M$ be the generic $g\times 2d$ matrix,
\[
M=\left(\begin{array}{cccc}
	x_{1 1}& x_{1 2}& \cdots  & x_{1, 2d} \\
	x_{2 1}& x_{2 2}& \cdots  & x_{1, 2d} \\
	\vdots &\vdots  &	  & \vdots    \\
	x_{g 1}& x_{g 2}& \cdots  & x_{g, 2d}   
      \end{array}\right)
\]
over the polynomial ring ${\mathbb C}[x_{ij}]$. We will
let $M_{(i_1,i_2,\ldots,i_g)},$  where $ i_1<i_2<\ldots<i_g,$
be the following submatrix of $M$.
\[
M_{(i_1,i_2,\ldots,i_g)}=\left(\begin{array}{cccc}
	x_{1,i_1}& x_{1,i_2}& \cdots  & x_{1,i_g} \\
	x_{2,i_1}& x_{2,i_2}& \cdots  & x_{1,i_g} \\
	\vdots &\vdots  &	  & \vdots    \\
	x_{g,i_1}& x_{g,i_2}& \cdots  & x_{g,i_g}
	\end{array}\right).
\]
Also let
\[
N=\left(\begin{array}{cccc}
	x_{1 1}& x_{1 2}& \cdots  & x_{1, g-1} \\
	x_{2 1}& x_{2 2}& \cdots  & x_{1, g-1} \\
	\vdots &\vdots  &	  & \vdots    \\
	x_{g 1}& x_{g 2}& \cdots  & x_{g, g-1}  
	\end{array}\right).
\]
Let $f$ be the product of the $(g-1)\times (g-1)$
minors of $N$. Let $R={\mathbb C}[x_{ij}]_f$. 
Let  $I$ be the ideal generated by the $g\times g$ minors
of $M$ in ${\mathbb C}[x_{ij}]$. Finally let $J$
be the ideal of ${\mathbb C}[x_{ij}]$ generated  by
the minors of the form ${\rm det}(M_{(1,2,\ldots,g-1,i)})$,
as $i$ ranges over, $g\le i\le 2d$. 
We wish to prove 

\begin{proposition} Consider the ideals $I_f,J_f$ obtained by extending
$I$ and $J$ to  the
ring $R$. We have $I_f=J_f$.
\label{P:reduce}
\end{proposition}

\begin{proof} It is clear that $J_f\subseteq I_f$. We proceed 
by showing that $\sqrt{J_f}=\sqrt{I_f}$ and then showing
that $J_f$ is equal to its radical.

We begin by showing $\sqrt{J\langle f\rangle}=
\sqrt{I\langle f\rangle}$. Here
$\langle f\rangle$ is the ideal generated by $f$. To show the above
it suffices to show that the two ideals have the same 
zero locus inside ${\mathbb A}^{g\times 2d}$. It is clear
that
\[ 
Z(J\langle f\rangle)=Z(J)\cup Z(f)\supseteq Z(I)\cup Z(f)=Z(I\langle 
f\rangle ).
\]
Now take $p=(p_{ij})$ in the zero locus of $J\langle f\rangle$. 
We may assume $p$ is not in the zero locus of $f$, for
otherwise we are done. Consider the matrix
\[
M_p=\left(\begin{array}{cccc}
		p_{11}& p_{12}&\ldots &p_{1,2d}\\
		p_{21}& p_{22}&\ldots &p_{2,2d}\\
		\vdots& \vdots&       &\vdots  \\
		p_{g1}& p_{g2}&\ldots &p_{g,2d}\\
	   \end{array}\right).
\]
Showing that $p\in Z(I)$ is equivalent to showing that 
${\rm rank}(M_p)\le g-1$. Since
$(p_{ij})\notin Z(f)$ the first $g-1$ columns of $M_p$
are linearly independent. As $(p_{ij})\in Z(J)$,
\[
{\rm det}((M_p))_{(1,2,\ldots,g-1,i)}=0,
\]
for $g\le i\le 2d$. So the $i$th column is
in the linear span of the first $g-1$ columns and
we are done. We have shown $\sqrt{J\langle f\rangle}=
\sqrt{I\langle f\rangle}$. An elementary argument now shows that
$\sqrt{I_f}=\sqrt{J_f}$.

Finally we need to show that $J_f$ is radical. Notice that
$J_f$ is generated by polynomials of the form
\[
\begin{array}{c}
{\rm det}(M_{(1,2,\ldots,g-1,i)})=
{\rm det}(N_1)x_{1i}-{\rm det}(N_2)x_{2i}+\ldots \\
(-1)^g{\rm det}(N_{g-1})x_{g-1,i}.
\end{array}
\]
Here $N_j$ is the submatrix of $N$ obtained by deleting
the $j$th row. Each of the ${\rm det}(N_j)$ are units
in our ring $R$. The result follows from the following lemma.
\end{proof}

\begin{lemma} 
Let $A$ be a reduced ring and consider the polynomial ring $B = A[x_{ij}]$,
where $1\le i\le n$ and $1\le j\le m$. Consider elements
\[
f_i = u_{i1}x_{i1} + u_{i2}x_{i2} + \ldots + u_{im}x_{im}.
\]
Form the ideal $I = (f_1,f_2,\ldots, f_n)$. If the $u_{ij}$
are units in $A$ then $B/I$ is reduced.
\end{lemma}
\begin{proof}
Observe that $B/I \cong A[x_{ij}]$ but with new index ranges $2\le i\le n$
and $2\le j\le m$.
\end{proof}

Now let $M$ be the matrix
\[
M = \left(\begin{array}{ccc}
	x_{11} & \cdots & x_{1,2d} \\
	\vdots &        & \vdots  \\
	x_{g1} & \cdots & x_{g,2d}
      \end{array}\right)
\]
over the polynomial ring ${\mathbb C}[x_{ij}]$.
Consider the submatrices
\[
A = \left(\begin{array}{ccc}
	x_{11} & \cdots & x_{1,d} \\
	\vdots &        & \vdots  \\
	x_{d,1} & \cdots & x_{d,d}
      \end{array}\right) \quad
B = \left(\begin{array}{ccc}
	x_{1,d+1}& \cdots & x_{1,2d} \\
	\vdots   & 	  & \vdots \\
	x_{d,d+1}& \cdots & x_{d,2d}
     \end{array}\right).
\]
Set  $f={\rm det}(A)$ and $g={\rm det}(B)$.
We will be interested in the following ideals of the ring 
${\mathbb C}[x_{ij}]_{fg}$. Let $I$ be ideal of the
$g\times g$ minors of $M$ and let $J$ be the ideal of the
$g\times g$ minors of 
\[
 N=   M \left(\begin{array}{cc}
	A^{-1} & 0 \\
	 0     & B^{-1}
	\end{array}\right).
\]

\begin{lemma}
\label{L:schemestructure}
The ideals $I$ and $J$ of ${\mathbb C}[x_{ij}]_{fg}$ are
equal.   
\end{lemma}
\begin{proof}
The subschemes of ${\rm spec}({\mathbb C}[x_{ij}]_{fg})$ 
defined by $I$ and $J$ are supported on the same
closed subset. So it suffices to show that both $I$
and $J$ are reduced. The fact that $I$ is reduced is the
fundamental theorem of invariant theory, see \cite{acgh}.
To show that $J$ is reduced consider the ${\mathbb C}$ algebra
automorphism of ${\mathbb C}[x_{ij}]_{fg}$ defined by
\[
x_{ij}\mapsto y_{ij}
\]
where
\[
M\left(\begin{array}{cc}
	A & 0 \\
	0 & B
	\end{array}\right)
=\left(\begin{array}{ccc}
	y_{11} & \cdots & y_{1,2d} \\
	\vdots &        & \vdots  \\
	y_{g1} & \cdots & y_{g,2d} 
	\end{array}\right).
\]
This automorphism carries $J$ to $I$ so we are done.
\end{proof}

\msection{A Subvariety of $\Grass{d-1,g-1}\times\Grass{d-1,g-1}$}

We let $\Grass{d-1,g-1}$ denote the Grassmanian paramaterizing
$(d-1)$ dimensional linear subspaces of ${\mathbb P}^{g-1}$. Let
\[
\begin{array}{c}
\f{d}{g}= \{(V,W)\in\Grass{d-1,g-1}\times\Grass{d-1,g-1}\mid \\
V\subseteq H, W\subseteq H\ 
{\rm for\ some\ hyperplane\ }H\subseteq
{\mathbb P}^{g-1}\}.
\end{array}
\]

In the above $V$ and $W$ are closed points of the Grassmanian.
We wish to describe the reduced scheme structure on $\f{d}{g}$. 
First we recall how to cover Grassmanian with open affines isomorphic
to ${\mathbb C}^{d(g-d)}$. 

Let $V\in\Grass{d-1,g-1}$ be a closed point. So $V$ can be thought of
 as the column space of a $g\times d$ matrix $A$. Write
\[
A=\left( \begin{array}{cccc}
		a_{11} & a_{12}& \cdots & a_{1d} \\
		a_{21} & a_{22}& \cdots & a_{2d} \\
           	\vdots & \vdots&  \vdots& \vdots \\
    		a_{g1} & a_{g2}& \cdots & a_{gd} \\
	\end{array} \right).
\] 
This representation is unique upto the action of 
${\rm GL}(d,{\mathbb C})$. 

 Let $I=( i_1,i_2,\ldots i_d)$, 
where $i_j\in\{1,2,\ldots ,g\}$ and   
$i_1<i_2<\ldots<i_d$. We will denote by $A^I$ the following
$d\times d$ submatrix of $A$:\footnote{In the preceeding 
section we defined $A_I$. In that section the submatrix $A_I$ 
of $A$ was obtained by choosing columns of $A$, while here
we are choosing rows.}
\[
A^I=\left( \begin{array}{cccc}
		a_{i_11} & a_{i_12}& \cdots & a_{i_1d} \\
		a_{i_21} & a_{i_22}& \cdots & a_{i_2d} \\
           	\vdots & \vdots&  \vdots& \vdots \\
    		a_{i_g1} & a_{i_g2}& \cdots & a_{i_gd} \\
	\end{array} \right).
\]

Now since the rank of $A$ is $d$, the matrix $A$ has a 
non vanishing $d\times d$ minor. Let this minor be
$\det(A^I)$. The matrix $A'=A(A^I)^{-1}$ also has column space
equal to $V$, furthermore it is the unique representative 
with $({A'})^I={\rm Id}_d$. For each $I=( i_1,i_2,\ldots i_d)$ as
above, set
\[
\begin{array}{c}
U_I=\{ V\in\grass{d}{g}\mid\ {\rm the\ }I{\rm\ minor\ of\ a\ 
matrix} \\  {\rm representative\ of\ }V{\rm\ is\ invertible}\}.
\end{array}
\]
There is a bijection $U_I\cong{\mathbb C}^{d.(g-d)}$, which
is in fact an isomorphism. For further details see \cite{griffithsharris}
or \cite{harris}. 

It follows from the above that $\Grass{d-1,g-1}\times\Grass{d-1,g-1}$
has an open affine cover constisting of opens of the form
 $U_I\times U_J\cong{\mathbb C}^{2d(g-d)}$.
Now take $(V,W)\in U_I\times U_J$, with $V$ the column space
of a matrix $A$ and $W$ the column space of a matrix $B$. The locus
we are trying to describe, $\f{d}{g}$, consists of those pairs
$(V,W)$ such that ${\rm rank}(A|B)<g$. Here $(A\mid B)$ is the matrix obtained
by augmenting the matrix $A$ with the matrix $B$.
Now the rank of  $(A|B)< g$ if and only if the $g\times g$
minors of $(A\mid B)$ vanish. The latter 
condition holds if and only if the $g\times g$ minors of the
matrix $(A \mid B)C$ vanish where
\[
C=\left(\begin{array}{cc} 
		A_I^{-1}& 0 \\
		0       & B_J^{-1} 
	\end{array}\right).
\] 
The entries of the matrix $(A B)C$ determine the image
of $(V,W)$ under the isomorphism 
$U_I\times U_J\cong {\mathbb C}^{d(g-d)+d(g-d)}.$
So the ideal generated by the $g\times g$ minors of
$(A\mid B)C$ determines a scheme structure on
$\f{d}{g}\cap U_I\times U_J$.It follows from \cite{acgh}
pg. 71 that this scheme structure is reduced, being a 
specialization of the ideal $I_k$ defined there. Hence
 these ideal sheaves on $U_I\times U_J$ glue together
to give an ideal sheaf for the reduced structure on $\f{d}{g}$.

We let 
\[
U_{\rm F}=\{ (V,W)\in\f{d}{g}\mid {\rm rank}(A\mid B)=g-1 \}.
\]
There is a morphism 
\[
\xymatrix{
\alpha\colon U_F \ar@{-->}[r]&\Pdual,}
\]
It takes a closed point $(V,W)$ to the linear span
of $V$ and $W$.
We will denote $\overline{U}_F$ by $\f{d}{g}_\m$.

\msection{ The construction of $\e{d}{g}$}         

In this section let $C$ be a  curve
of genus $g\ge 4$. Let $(g-1)/2<d<g-1$. We have a morphism
\[
\Gauss\times\Gauss\colon\w_\smooth\times\w_\smooth\rightarrow
\Grass{d-1,g-1}\times\Grass{d-1,g-1}.
\]
Define $\e{d}{g}\hookrightarrow\w_\smooth\times\w_\smooth$
to be the fibre over  $\f{d}{g}$. We take $U_{\rm E}$ to be
the preimage of $U_{\rm F}$ and $\e{d}{g}_\m$ to be the closure
of $U_{\rm E}$.There is  morphism
\[
\beta\colon U_{\rm E}\longrightarrow ({\mathbb P}^{g-1})^*.
\]
We have, by theorem \ref{T:Gaussmap},
\[
\beta(u(D),u(D'))=\overline{{\phi_K(D)}\cup{\phi_K(D')}},
\]
where $(D,D')\in\symm\times\symm$ are divisors whose image 
under the Abel-Jacobi map is in $\w_\smooth$. Recall that
$\overline{A}$ means linear span of some 
subset $A$ of ${\mathbb P}^{g-1}$ in ${\mathbb P}^{g-1}$.
Notice that $\lintd{D}{D'}$ is a hyperplane in ${\mathbb P}^{g-1}$,
for the condition $\abt{D}{D'}\in\e{d}{g}$ forces $\lintd{D}{D'}$
to be contained in a hyperplane and the condition $\abt{D}{D'}\in U_{\rm E}$
forces $\lintd{D}{D'}$ to be exactly a hyperplane.  

A generic hyperplane $H\in\Pdual$ intersects $C$ in $2g-2$ points
that are in general position , see \cite{acgh}. So suppose that 
$H.C=p_1+p_2+\ldots+p_{2g-2}$. Then by \ref{T:geometricRR}, the
pair 
\[
\abt{p_1+p_2+\ldots+p_d}{p_{d+1}+p_{d+2}+\ldots+p_{2d}},
\]
(notice $2d<2g-2$) is a closed point of $\w_\smooth\times\w_\smooth$. 
Furthermore the above pair, gives a point in $U_E$ mapping
to $H$ under $\beta$. Hence $\beta$ is dominant. Since a
hyperplane can only intersect $C$ in a finite number of points,
the map $\beta$ is quasi-finite. It follows that $U_E$ has
dimension $g-1$. 

We let $C^*$ denote the dual variety to $\phi_K(C)$.

\begin{lemma}
\label{L:esmooth} 

\noindent (a) Suppose that $C$ is a non-hyperelliptic
curve. Let $H\in\Pdual-C^*$. If 
$\beta(\abt{D}{D'})=H$ 
then $\abt{D}{D'}$ lies on a component of $\e{d}{g}$ of 
dimension $g-1$ and is in the
smooth locus of $\e{d}{g}_\m$.

\noindent (b) Suppose that $C$ is hyperelliptic. Let 
$H\in\Pdual-C^*$ and assume also that $H$ does not pass
through any of the branch points of the canonical map
$\phi_K\colon C\rightarrow{\mathbb P}^{g-1}$. If 
$\beta(\abt{D}{D'})=H$ 
then $\abt{D}{D'}$ lies on a component of $\e{d}{g}$ of
dimension $g-1$ and is in the
smooth locus of $\e{d}{g}_\m$.
\end{lemma}

\begin{proof}

The following proof is for (a).

Write $D=p_1+p_2+\ldots+p_d$ and $D'=p'_1+p'_2+\ldots+p'_d$. 
We choose local coordinates $z_i$ and $z'_i$ on $C$ centred
at $p_i$ and $p'_i$ respectively. Now as $H$ is not a tangent
hyperplane $C$, we have $p_i\ne p_j$ and $p'_i\ne p'_j$ for
$i\ne j$. It follows that $z_1, z_2\ldots, z_d$ and
 $z'_1,z'_2\ldots  , z'_d$ descend
to local co-ordinates on $\symm\times\symm$ centred at $(D,D')$.
Furthermore, by \ref{T:abels}, the Abel-Jacobi map is an
isomorphism around $(D,D')$, since $u(D),u(D')\in\w_\smooth$.
 So we have some local co-ordinates
on $\w\times\w$ centred at $\abt{D}{D'}$. Let $\omega_1,\ldots
\omega_g$ be a basis for ${\rm H}^0(\Omega^1_C)$. We write
$\omega_j$ as $\Omega_{ji}(z_i)dz_i$ in a neighbourhood of 
$p_i$ and as $\Omega'_{ji}(z'_j)dz'_j$ in a neighbourhood of
$p'_j$. Let
\[ 
M(z)= \left(\begin{array}{cccccc}
\Omega_{11}(z_1) & \ldots &\Omega_{1d}(z_d) &\Omega'_{11}(z'_1)&\ldots
&\Omega'_{1d}(z'_d) \\
\Omega_{21}(z_1) & \ldots &\Omega_{2d}(z_d) &\Omega'_{21}(z'_2)&\ldots
&\Omega'_{2d}(z'_1) \\
\vdots& &\vdots&\vdots& &\vdots \\
\Omega_{g1}(z_1) & \ldots &\Omega_{gd}(z_d) &\Omega'_{g1}(z'_1)&\ldots
&\Omega'_{gd}(z'_d) \\
	\end{array}\right).
\]
In a neighbourhood of $\abt{D}{D'}$, $\e{d}{g}$ is defined by the
vanishing of the $g\times g$ minors of $M(z)$, by \ref{L:schemestructure}.
Now by \ref{T:geometricRR}, ${\rm dim}\lind{D}=d-1$, so in a neighbourhood
of $\abt{D}{D'}$ the first $d$ columns of $M(z)$ are linearly
independent. Since $M(Z)$ has rank $g-1$ at the point $\abt{D}{D'}$
we may reindex the points of $D'$ so that the first $g-1$ columns
of $M(z)$ are linearly independent in a neighbourhood of 
$\abt{D}{D'}$. Set
\[
f_i={\rm det}M(z)_{(1,2,\ldots,g-1,i)},
\]
where $g-1<i\le 2d$. By \ref{P:reduce}, $\e{d}{g}$ is 
defined by $f_i$ in a neighbourhood of $\abt{D}{D'}$. The assertion
that $\abt{D}{D'}$ lies on a component of dimension $g-1$ of
$\e{d}{g}$ follows.

By definition, $f_j$ is independent of the co-ordinates $z'_i$ for
$g-d\le i\le d$ and $i\ne j-d$. So the Jacobian matrix is of the form
\[
\left.\left(\begin{array} {cccc}
\vspace{1pt}
\pde{f_g}{z_1} & \pde{f_{g+1}}{z_1} & \cdots  & \pde{f_{2d}}{z_1} \\ 
\vspace{1pt}
\pde{f_g}{z_2} & \pde{f_{g+1}}{z_2} & \cdots  & \pde{f_{2d}}{z_2} \\
\vdots         &  \vdots            &         & \vdots \\
\vspace{1pt}
\pde{f_g}{z_d} & \pde{f_{g+1}}{z_d} & \cdots  & \pde{f_{2d}}{z_d} \\
\vspace{1pt}
\pde{f_g}{z'_1} & \pde{f_{g+1}}{z'_1} & \cdots  & \pde{f_{2d}}{z'_1} \\ 
\vspace{1pt}
\pde{f_g}{z'_2} & \pde{f_{g+1}}{z'_2} & \cdots  & \pde{f_{2d}}{z'_2} \\
\vdots         &  \vdots            &         & \vdots \\
\vspace{1pt}
\pde{f_g}{z'_{g-d-1}} & \pde{f_{g+1}}{z'_{g-d-1}} & \cdots  & \pde{f_{2d}}{z'_{g-d-1}} \\
\vspace{1pt}
\pde{f_g}{z'_{g-d}} & 0 & \cdots  & 0 \\
\vspace{1pt}
 0 & \pde{f_{g+1}}{z'_{g-d+1}} & \cdots  & 0 
\\
\vspace{0pt}
\vdots		&	\vdots	      &		& \vdots	\\
\vspace{0pt}
 0 & 0& \cdots  & \pde{f_{2d}}{z'_d} 
\end{array}\right)\right|_{\abt{D}{D'}}.
\]
Suppose that $\abt{D}{D'}$ is a singular point of ${\rm E}(d,g)$. This is true
 if and only if
the above matrix has rank smaller than $2d-g+1$. 
It has rank smaller than $2d-g+1$ if and only if 
\[
\left.\pde{f_j}{z'_j}\right|_{\abt{D}{D'}}=0
\]
 for some $j$.
Now

\begin{eqnarray}
0& = &\left.\pde{f_j}{z'_j}\right|_{\abt{D}{D'}} \nonumber \\ \nonumber 
 & = &
\left|\begin{array}{ccccccc}
\vspace{5pt}
\Omega_{11}(p_1) \hspace{-5pt} &\cdots \hspace{-5pt} & \Omega_{1d}(p_d) \hspace{-5pt} &
\Omega'_{11}(p'_1) \hspace{-5pt} &\cdots \hspace{-5pt} &
\Omega'_{1,g-1-d}(p'_{g-1-d})\hspace{-5pt} &
\left.\pde{\Omega'_{1j}}{z'_j}\right|_{p_j} \\
\vspace{5pt}
\Omega_{21}(p_1) \hspace{-5pt}&\cdots\hspace{-5pt} & \Omega_{2d}(p_d)\hspace{-5pt} &
\Omega'_{21}(p'_1)\hspace{-5pt} &\cdots\hspace{-5pt} &
\Omega'_{2,g-1-d}(p'_{g-1-d})\hspace{-5pt} &
\left.\pde{\Omega'_{2j}}{z'_j}\right|_{p_j} \\
\vspace{5pt}
\vdots\hspace{-5pt} &\hspace{-5pt} & \vdots\hspace{-5pt} & \vdots\hspace{-5pt} & \hspace{-5pt}& \vdots\hspace{-5pt} & \vdots \\
\vspace{5pt}
\Omega_{g1}(p_1)\hspace{-5pt} &\cdots\hspace{-5pt} & \Omega_{gd}(p_d)\hspace{-5pt} &
\Omega'_{g1}(p'_1)\hspace{-5pt} &\cdots\hspace{-5pt} &
\Omega'_{g,g-1-d}(p'_{g-1-d})\hspace{-5pt} &
\left.\pde{\Omega'_{gj}}{z'_j}\right|_{p_j}
\end{array}\right|
\end{eqnarray}

The first $g-1$ columns lie inside $H$.
So it follows that the last column is contained in $H$. 
This implies the tangent line to $p'_j$ is in $H$, which
in turn contradicts $H\not\in C^*$.

A similar argument proves (b).
\end{proof}

\msection{Generic Tangent Hyperplanes}

Let $C$ be a curve with a fixed non-degenerate embedding
$\phi\colon C\hookrightarrow{\mathbb P}^n$, with $n\ge 3$.
Recall that all curves are assumed to be smooth
and projective. The genus of our curve will also be assume to
be $\ge 4$. We will denote by $C^*$ the 
dual variety to $C$ inside $({\mathbb P}^n)^*$. By
forming the incidence correspondence
\[
\Sigma=\{(p,H)\mid p\in C,\ H\in({\mathbb P}^n)^*,\ 
{\rm T}_p(C)\subseteq H\}
\]
and using standard arguments we see that $C^*$ is an
irreducible hypersurface in $({\mathbb P}^n)^*$. We use the
notation ${\rm T}_p(C)$ to denote the tangent line 
to $C$ at $p$ inside ${\mathbb P}^n$. 

Let $\phi_2\colon C\rightarrow\Grass{2,n}$ be the
second associated curve to $\phi$. So $\phi_2(p)$
is the unique plane having intersection order at
least $3$ with $C$ at $p$. (See \cite{griffithsharris}, pg. 263).
Let $\Gamma_2\subseteq C\times\Grass{2,n}$ be the 
graph of $\phi_2$.
We form the incidence correspondence 
\[
\Sigma''=\{(p,P,H)\in\Gamma_2\times({\mathbb P}^n)^*
\mid (p,P)\in\Gamma_2,\ {\rm and}\ P\subseteq H\}.
\]
Let $\pi_C\colon$ be the projection $\pi_C\colon
\Sigma''\rightarrow C$. The fibre over $p\in C$ is
irreducible of dimension $n-3$. It follows that 
$\Sigma''$ is irreducible of dimension $n-2$. The
projection from $\Sigma''$ to $C^*$ is a finite
morphism, hence the locus of hyperplanes having 
intersection at least $3$ at some point of $C$ is
an irreducible closed subvariety of codimension $1$
inside $C^*$. 

\begin{lemma}
\label{L:gentangent}
Let $\phi_K\colon C\rightarrow{\mathbb P}^{g-1}$
be the canonical morphism. 

\noindent (a) Suppose that $C$ is a non-hyperelliptic
curve so that $\phi_K$ is an immersion. Then for
a generic $H\in C^*\subseteq\Pdual$,
\[
H.C=2p_1+p_2+p_3+\ldots+p_{2g-3}
\]
where the $p_i$ are distinct.

\noindent (b) Suppose that $C$ is a hyperelliptic curve 
so that $\phi_K(C)$ is a rational normal curve . Let $C^*$
be the dual variety to $\phi_K(C)$.
Let $b_1,\ldots,b_{2g+2}$ be the
branch points of $\phi_K$. We denote by $b_i^*\subseteq\Pdual$
the dual variety to $b_i$, consisting of all hyperplanes
through $b_i$. So $b_i^*$ is a hyperplane in $\Pdual$. 
Then for a generic 
\[
H\in C^*\cup b_1^*\cup\ldots\cup b_{2g+2}^*
\]
we have that
\[
H.C=2p_1+p_2+p_3+\ldots+p_{2g-3}
\]
where the $p_i$ are distinct.
\end{lemma}
\begin{proof}

\noindent(a)  We have seen, in the discussion preceeding the
lemma, that for a generic $H\in C^*$,
$H.C$ has no points of multiplicity $3$.
So we need to show that a generic tangent hyperplane
has only one point of multiplicity 2. Form the
incidence correspondence
\[
\begin{array}{c}
\Sigma'=\{(p,q,P,H)\in C\times C\times\Grass{3,g-1}\times\Pdual
\mid p\ne q\   \\ {\rm T}_p(C),{\rm T}_q(C)\subseteq P\subseteq H \}.
\end{array}
\]
Note that $\Sigma'$ is only locally closed in $C\times C\times
\Grass{3,g-1}\times\Pdual$. Let
\[
\Sigma'=\Sigma'_1\cup\Sigma'_2\cup\ldots\cup\Sigma'_l
\]
be an irreducible decomposition for $\Sigma'$. There
is a projection $\Sigma'\rightarrow C\times C$. 
From \cite{hartshorne} IV Theorem (3.10) there is a 
closed subset $X\subseteq C\times C$ such that
for each $(p,q)\not\in D$, the tangent lines
${\rm T}_p(C)$ and ${\rm T}_q(C)$ do not meet and $X$
has codimension $1$ in $C\times C$. Consider
the restricted projection
\[
\pi_i\colon\Sigma'_i\rightarrow C\times C.
\]
Now if there is a point $(p,q)\not\in X$, and
in the image of $\Sigma'_i$, the fibre over $(p,q)$
has dimension $g-5$ as ${\rm T}_p(C)$ and ${\rm T}_q(C)$
span a 3-plane in ${\mathbb P}^{g-1}$. Hence
\[
\begin{array}{ccc}
{\rm dim}\hspace{3pt} \Sigma'_i&\le & {\rm dim}\hspace{3pt}
C\times C + {\rm dim}(
{\rm fibre}) \\
 &= & g-3.
\end{array}
\]
(Note that if $g=4$, then there is no such $(p,q)$.) 
If there is no such $(p,q)$ then the projection can be 
factored as
\[
\pi_i\colon\Sigma'_i\rightarrow X.
\]
Now the fibre over a point has dimension $g-4$. So
as above ${\rm dim}\Sigma'_i\le g-3$. Hence, for the
closure $\overline{\Sigma'}$, we have 
\[
{\rm dim}\overline{\Sigma'}\le g-3.
\]
So the image of the projection 
$\overline{\Sigma'}\rightarrow C^*$ has smaller than
dimension $g-2$. Since $C^*$ is a hypersurface,
the result follows.

\noindent(b) First consider $H\in C^*$. 

By the remark proceeding the lemma, it suffices to show
that for a generic $H\in C^*$, $H.C$ has
only one point of multiplicity two and $H$ does
not pass through one of the $b_i$. The first assertion 
follows as in (a). For the second assertion notice that
$C^*,b_1^*,\ldots,b_{2g+2}^*$ are distinct hypersurfaces
in $\Pdual$. The result follows. 

This last remark also deals with the case $H\in {b_i}^*$.
\end{proof}

\msection{Proof of the Generalized Torelli Theorem}

We wish to prove

\begin{theorem}

Let $C$ be a smooth projective curve over ${\mathbb C}$
of genus $g\ge 1$. If $1\le d\le g-1$ is an integer then
the pair $(J(C),\w)$ determine the curve, that is if
$(J(C),\w(C))\cong(J(C'),\w(C'))$ for some other
smooth projective curve $C'$ then $C'\cong C$. 
\end{theorem}

\begin{proof} We may assume $g\ge 4$ as the cases $g=1,2,3$ are
covered by the regular Torelli theorem. Furthermore we may 
reduce to the case $(g-1)/2< d<g-1$ as follows. If $d=g-1 $ we
are done by Torelli's theorem. If $d< (g-1)/2$ then
choose $n$ so that $(g-1)/2<nd\le g-1$. Now we may replace
$\w$ by
\[
W^{nd}=\underbrace{W^d+W^d+\ldots +W^d.}_{n\ {\rm times}}
\]

We will study the branch locus of the map 
\[ \xymatrix{
\beta\colon\e{d}{g}_{{\rm main}} \ar@{-->}[r] & \Pdual.}
\]

 Note that we can recover the rational map
$\beta$ from the information $(J(C),\w)$. Now let
$U_{\rm E}\subseteq\e{d}{g}$ be the open subset defined at the
start of \S 5. We have a morphism
$\beta|_{U_{\rm E}}\colon U_{\rm E}\rightarrow\Pdual$. Let $B$ be
the branch locus of $\beta$. This is the image of the ramification
locus inside $\Pdual$. A closed point $p$ is in the ramification
locus if and only if $\beta$ fails to be a local analytic 
isomorphism at $p$. At this point we break the proof into two
cases, the case where $C$ is non-hyperelliptic and the case where
$C$ is hyperelliptic. 

First we study the case where $C$ is non-hyperelliptic.
 We will show that $\bar{B}=C^*$. Then
$C$ can be recovered from this information, see \cite{harris}.

First we show that $\bar{B}\subseteq C^*$. Let $H\not\in C^*$. Then 
$H.C=p_1+p_2+\ldots +p_{2g-2}$ with the $p_i$ distinct.
Let $T\subseteq\Pdual$ be all the hyperplanes having transverse
intersection with $C$, that is $T=\Pdual - C^*$. The
incidence correspondece
\[
I=\{(p,H)\in C\times T\mid p\in{\rm Supp\ }H.C \}\rightarrow T
\]
is a $(2g-2)$-sheeted covering space of $T$,\cite{acgh} pg.110. Given
$\abt{D}{D'}\in U_E$ with $\beta(\abt{D}{D'})=H\in T$. It is
claimed that there exists an open neighbourhood $V$ in
the usual topology such that
\[
\beta|_V\colon V\rightarrow\beta(V)
\]
is an injection. To see this, first take
$H\in W\subseteq T$, with sheets $W_1,W_2,\ldots , W_{2g-2}$. Let
$\mu_i$ be the compostion $W\rightarrow W_i\rightarrow C$,
which is holomorphic. Write $D=p_1+\ldots +p_d$. The $p_i$
are distinct by choice of $H$, so we may find opens $p_i\in U_i$
such that

(1) $U_i\cap U_j = \empty$ for $i\ne j$

(2) $U_i\subseteq\mu_j(W)$ for some $j$.

Writing $D'=p'_1+\ldots +p'_d$ we may find similar
opens $U'_i$. Set $U=U_1\times\ldots\times U_d$, 
$U'=U'_1\times\ldots\times U'_d$. By condition (1), 
$U\times U'$ is an open neighbourhood of $(p_1+\ldots +p_d,p'_1+\ldots
+p'_d)$ on $C^{(d)}\times C^{(d)}$. As the Abel-Jaacobi map
is an isomorphism near $(p_1+\ldots +p_d,p'_1+\ldots
+p'_d)$, as $\abt{D}{D'}\in\w{d}_\smooth\times\w{d}_\smooth$.
We take $V=\beta^{-1}\cap (U\times U')\cap U_E$. It is
easy to see that this works.

 It follows from theorem 7.6, of \cite{grauert},
that $\beta$ is a local isomorphism at $\abt{D}{D'}$ since
this point is in the smooth locus of $\e{d}{g}$ by lemma 
\ref{L:esmooth}. It remains to show that $B$ contains an open
dense subset of $C^*$.

By \ref{L:gentangent} there exists an open subset $V\subseteq C^*$
 such that for each $H\in V$,
\[
H.C=2p_1+p_2+\ldots +p_{2g-3},
\]
with the $p_1,\ldots,p_{2g-3}$ are distinct. Since $g\ne 0$ and
\[
K\sim 2p_1+p_2+\ldots +p_{2g-3},
\]
we have that $H=\overline{\phi_K(p_1+p_2+\ldots +p_{2g-3})}$. (Notice that
there is no $2$ in front of $p_1$ in the last statement.)
After reindexing we may assume that $p_1,p_2,\ldots,p_{g-1}$
span $H$ and the tangent line at $p_1$ to $C$ lies inside
$H$. Let 
\[
D=q_1+q_2+\ldots+q_{d}\quad{\rm and}\quad D'=q'_1+q'_2+
\ldots+q'_d
\]
where $q_i=p_i$ for $1\le i\le d$ and $q'_i=p_{g-i}$ for
$1\le i\le d$. So $\abt{D}{D'}\in U_{\rm E}$. Let $z_i$ (resp. $z'_i$)
be local coordinates centred at $q_i$ (resp. $q'_i$). 
Since $q_i\ne q_j$ (resp. $q'_i\ne q'_j$) for $i\ne j$, we
have local coordinates $(z_1, z_2, \ldots ,z_d)$
(resp. $(z'_1, z'_2,\ldots ,z'_d$) on $\symm$ centred
at $(q_1,q_2,\ldots,q_d)$ (resp. $(q'_1,\ldots,q'_d)$). As
$u$ is an isomorphism around $D$ (resp. $D'$), by \ref{T:abels} 
and \ref{T:geometricRR} and as ${\rm dim}\overline{\phi_K(D)}=d-1$
(resp. ${\rm dim}\overline{\phi_K(D')}=d-1$), we have that 
$((z_1, z_2,\ldots z_d),(z'_1, z'_2,\ldots z'_d))$
descend to local coordinates on $\w\times\w$ centred at
$\abt{D}{D'}$.

Choose a basis $\omega_1,\ldots,\omega_g$ for
${\rm H}^0(C,\Omega^1_C)$ and write $\omega_i=
\Omega_{ij}(z_j)dz_j$ (resp. $\omega_i=
\Omega'_{ij}(z'_j)dz'_j$). Let

\[
\begin{array}{l}
M(z)= \\ \left(\begin{array}{cccccccc}
\Omega_{11}(z_1)\hspace{-5pt} & \Omega_{12}(z_2)\hspace{-5pt}  &\ldots\hspace{-5pt}  &\Omega_{1,d}(z_{d})\hspace{-5pt} &
\Omega'_{11}(z'_1)\hspace{-5pt}  & \Omega'_{12}(z'_2)\hspace{-5pt}  &\ldots\hspace{-5pt} 
&\Omega'_{1,d}(z'_{d})\\
\Omega_{21}(z_1)\hspace{-5pt}  & \Omega_{22}(z_2)\hspace{-5pt}  &\ldots\hspace{-5pt}  &\Omega_{2,d}(z_{d})\hspace{-5pt} &
\Omega'_{21}(z'_1)\hspace{-5pt}  & \Omega'_{22}(z'_2)\hspace{-5pt}  &\ldots\hspace{-5pt} 
&\Omega'_{2,d}(z'_{2})\\
\ldots \hspace{-5pt} & \ldots\hspace{-5pt}  & \hspace{-5pt}  &\ldots\hspace{-5pt}  &\ldots \hspace{-5pt} 
& \ldots \hspace{-5pt}&\hspace{-5pt} & \ldots \\
\Omega_{g1}(z_1)\hspace{-5pt} & \Omega_{g2}(z_2)\hspace{-5pt} &\ldots\hspace{-5pt} &\Omega_{g,d}(z_{d})\hspace{-5pt}&
\Omega'_{g1}(z'_1)\hspace{-5pt} & \Omega'_{g2}(z'_2)\hspace{-5pt} &\ldots\hspace{-5pt}
&\Omega'_{g,d}(z'_{d})
	\end{array}\right)
\end{array}
\]
and let 
\[
\begin{array}{l}
M'(z)= \\ \left(\begin{array}{cccccccc}
\pde{\Omega_{11}(z_1)}{z_1}\hspace{-5pt} & 
\Omega_{12}(z_2)\hspace{-5pt} &\ldots\hspace{-5pt} &\Omega_{1,d}(z_{d})\hspace{-5pt}&
\Omega'_{11}(z'_1)\hspace{-5pt} & \Omega'_{12}(z'_2)\hspace{-5pt} &\ldots\hspace{-5pt}
&{\Omega'_{1,d}(z'_{d})}{z_1}\\
\pde{\Omega_{21}(z_1)}{z_1}\hspace{-5pt} & 
\Omega_{22}(z_2)\hspace{-5pt} &\ldots \hspace{-5pt}&\Omega_{2,d}(z_{d})\hspace{-5pt}&
\Omega'_{21}(z'_1)\hspace{-5pt} & \Omega'_{22}(z'_2)\hspace{-5pt} &\ldots\hspace{-5pt}
&\Omega'_{2,d}(z'_{2})\\
\ldots\hspace{-5pt} & \ldots\hspace{-5pt} & \hspace{-5pt} &\ldots\hspace{-5pt} &\ldots\hspace{-5pt}
 & \ldots\hspace{-5pt} & \hspace{-5pt}& \ldots \\
\pde{\Omega_{g1}(z_1)}{z_1}\hspace{-5pt} & 
\Omega_{g2}(z_2)\hspace{-5pt} &\ldots\hspace{-5pt} &\Omega_{g,d}(z_{d})\hspace{-5pt}&
\Omega'_{g1}(z'_1)\hspace{-5pt} & \Omega'_{g2}(z'_2)\hspace{-5pt} &\ldots\hspace{-5pt}
&\Omega'_{g,d}(z'_{d})
	\end{array}\right)
\end{array}
\]

By definition of $D$ and $D'$ the first $g-1$ columns of $M(z)$
are linearly independent in a neighbourhood of $\abt{D}{D'}$.
So $\e{d}{g}$ is defined by 
\[
f_i={\rm det}(M(z)_{1,2,\ldots,g-1,i}),
\]
where $g\le i\le 2d$, in a neighbourhood of $\abt{D}{D'}$.
(To see this, use \ref{P:reduce} as in \ref{L:esmooth})
Now since the tangent line to $C$ at $p_1$ is inside
$H$ we have 
\[
\pde{f_i}{z_1}={\rm det}(M'(z)_{1,2,\ldots,g-1,i})|_\abt{D}{D'}
=0.
\]
So the Jacobian matrix, as in the proof of \ref{L:esmooth},
reduces to 
\[
\left(\begin{array}{cccc}
0 & 0 & \cdots & 0 \\
\pde{f_g}{z_2} & \pde{f_{g+1}}{z_2} & \cdots &
\pde{f_{2d}}{z_2} \\
\vdots & \vdots & &\vdots \\
\pde{f_g}{z_d} & \pde{f_{g+1}}{z_d} & \cdots &
\pde{f_{2d}}{z_d} \\
\pde{f_g}{z'_1} & \pde{f_{g+1}}{z'_1} & \cdots &
\pde{f_{2d}}{z'_1} \\
\vdots & \vdots & &\vdots \\
\pde{f_g}{z'_{g-1-d}} & \pde{f_{g+1}}{z'_{g-1-d}} & \cdots &
\pde{f_{2d}}{z'_{g-1-d}} \\
\pde{f_g}{z'_{g-d}} & 0 & \cdots & 0 \\
0 & \pde{f_{g+1}}{z'_{g-d+1}} & \cdots & 0 \\
 & & \ddots & \\ 
0 & 0 & \cdots & \pde{f_2d}{z'_d} \end{array}\right)
\]

Arguing as in \ref{L:esmooth} we find that $\abt{D}{D'}$ is
a smooth point of $\e{d}{g}$. We also see that 
$\pde{ }{z_1}|_{\abt{D}{D'}}$ is
in the null space of the above Jacobian. Hence $\pde{ }{z_1}|_{\abt{D}{D'}}$ 
is in fact a tangential to $\e{d}{g}$ at $\abt{D}{D'}$. In order to show that 
$H\in B$ it will suffice to show that $\pde{ }{z_1}|_{\abt{D}{D'}}$
maps to zero under the morphism of tangent space induced by $\beta$.
 Let
\[
N(z)=\left(\begin{array}{cccccc}
\Omega_{11}(z_1) & \cdots & \Omega_{1d}(z_d) & 
\Omega'_{11}(z'_1) & \cdots & \Omega_{1,g-1-d}(z_{g-1-d}) \\
\vdots & & \vdots & \vdots & & \vdots \\
\Omega_{g1}(z_1) & \cdots & \Omega_{gd}(z_d) & 
\Omega'_{g1}(z'_1) & \cdots & \Omega_{g,g-1-d}(z_{g-1-d}) 
\end{array}\right).
\]

So $N(z)$ is just the first $g-1$ columns of $M(z)$. In
a neighbourhood of $\abt{D}{D'}$ the morphism 
$\beta\colon U\rightarrow\Pdual $ is given by 
$z\mapsto{\rm col.\ space}N(z)$. Identify $\Pdual\cong
{\mathbb P}(\bigwedge^{g-1}{\mathbb C}^g)$ we see that 
$\beta$ is the morphism 
\[
z\mapsto[{\rm det}(N(z)_1):{\rm det}(N(z)_2):\ldots:
{\rm det}(N(z)_g)].
\]
Recall that $N(z)_i$ is the submatrix of $N(z)$ obtained by 
deleting the $i$th row. We may assume that ${\rm det}(N(z)_1)
\ne 0$. So we need to show that 
\[
\pde{ }{z_1}|_{\abt{D}{D'}}(\frac{{\rm det}(N(z)_i)}{{\rm
det}(N(z)_1)}=0.
\]
That is 
\[
\pde{{\rm det}(N(z)_1)}{z_1}.{\rm det}
(N(z)_i)=\pde{{\rm det}(N(z)_i}{z_1}.
{\rm det}(N(z)_1)
\]
after evaluation at $\abt{D}{D'}$.
Let
\[
\pde{N(z)}{z_1}
\]
be the matrix obtained from $N(z)$ by differentiating the first
column with respect $z_1$.
Observe that 
\[
{\rm col.\ space}\pde{N(z)}{z_1}|_{\abt{D}{D'}}\subseteq
{\rm col.\ space} N(z)|_{\abt{D}{D'}}
\]
as the tangent line at $p_1$ lies inside $H$. It is a general
fact from linear algebra that given two $g\times (g-1)$ 
matrices $M,N$ with ${\rm col.\ space}M\subseteq{\rm col.\ space}
N$ then for each $i$ $j$ in the range $1\le i,j\le g$ we
have
\[
{\rm det}(M_i){\rm det}(N_j)={\rm det}(M_j){\rm det}( N_i).
\]
We will include the proof of this statement at the 
end of this proof for completeness. This shows
that $\bar{B}=C^*$.

Now we treat the case that $C$ is a 
hyperelliptic curve. We show that $\bar{B}=C^*\cup
b_1^*\cup b_2^*\cup\ldots\cup b_{2g+2}^*$ where the
$b_i$ are the branch points of the canonical 
morphism $\phi_K\colon C\rightarrow\Pdual$. The
proof is almost identical to the above. Here are a few
details. The same proof as in the non-hyperelliptic 
case shows that $\bar{B}\subseteq C^*\cup
b_1^*\cup b_2^*\cup\ldots\cup b_{2g+2}^*$, and similarly
we show that $\bar{B}\supseteq C^*$. To show that
$\bar{B}\supseteq{b_i^*}$ proceed as follows. From \ref{L:gentangent}
we know that for a generic $H\in b_i^*$ that
\[
H.C=2p_1 +p_2+\ldots+p_{2g-3}
\]
where the $p_i$ are distinct and $p_1=b_1$. As above we form, after
appropriate reindexing, 
\[
D=q_1+q_2+\ldots+q_d\quad{\rm and}\quad
D'=q'_1+q'_2+\ldots+q'_d.
\]
Note, these two divisors are defined exactly as they were
before. Also define, as before, $z_i$, $z'_i$, $M(z)$, $M'(z)$ 
and $f_i$. To see that 
\[
\pde{f_i}{z_1}|_{\abt{D}{D'}}=0,
\]
first observe that since $p_1$ is a branch point,
$J(\phi_K)|_{q_1}=0$. Around $q_1$, 
\[
\phi_K=[\Omega_{11}(z_1):\ldots:\Omega_{g1}(z_1)].
\]
We may assume that $\Omega_{11}(z_1)\ne 0$. Since the 
Jacobian at $q_1$ vanishes we see that
\[
\Omega_{11}(q_1)\pde{\Omega_{1j}(q_1)}{z_1}|_{q_1}=
\pde{\Omega_{11}(z_1)}{z_1}|_{q_1}\Omega_{1j}(q_1),
\]
which in turn implies
\[
\phi_K(q_1)=[\Omega_{11}(q_1):\ldots:\Omega_{g1}(q_1)]=
[\pde{\Omega_{11}}{z_1}:\ldots:\pde{\Omega_{g1}}{z_1}]\mid_{q_1}.
\]
So
\[
\pde{f_i}{z_1}|_{\abt{D}{D'}}=f_i|_{\abt{D}{D'}}=0.
\]
Now proceed as before.

\end{proof}

Here is the linear algebra result that was
needed before.

\begin{lemma} Let $M$, $N$ be two $g\times (g-1)$ 
matrices over ${\mathbb C}$. If 
\[
{\rm col. space} M\subseteq{\rm col. space N}
\]
then 
\begin{equation}
\label{E:det}
{\rm det}M_i{\rm det}N_j={\rm det}M_j{\rm det}N_i,
\end{equation}
for each $i,j$ with $1\le i,j\le g$. Recall that
$M_i$ is the submatrix of $M$ obtained by deleting
the $i$th row.
\end{lemma}

\begin{proof}

Firstly if ${\rm rank}M<g-1$ then both sides of 
(\ref{E:det}) vanish. So we may assume $M$, $N$ 
are of maximal rank and that there column spaces 
are equal. So $N=M.H$ for some $H\in{\rm Gl}(g-1,{\mathbb C})$.
The result follows from the observation $(M.H)_i=M_i.H$.
\end{proof}

\begin{corollary}

Let $C$ and $C'$ be two smooth projective curves and let 
$d$ be an integer less than or equal to the genus of $C$.
If $C^{(d)}\cong {C'}^{(d)}$ then $C\cong C'$.
\end{corollary}

\begin{proof} This is because the Albanese varaiety
${\rm Alb}(C^{(d)})$ is isomorphic to $J(C)$ and the image
of $C^{(d)}$ under the Albanese map is $W^d$.
\end{proof}

\end{document}